\documentclass[11pt]{article}%
\usepackage{graphicx}
\usepackage{amsmath}
\usepackage{amsfonts}
\usepackage{amssymb}
\providecommand{\U}[1]{\protect\rule{.1in}{.1in}}
\newtheorem{theorem}{Theorem}

\newtheorem{corollary}[theorem]{Corollary}

\newtheorem{lemma}[theorem]{Lemma}

\newtheorem{proposition}[theorem]{Proposition}

\newcommand{\dem}{\noindent \mbox{\bf Proof. }}
\newcommand{\QED}{\hspace*{\fill}\mbox{\small \bf  Q.E.D.}\bigskip}

\newcommand{\Nat}{{\mathrm{I\hspace{-.15em}N}}}

\renewcommand{\leq}{\leqslant}
\renewcommand{\geq}{\geqslant}

\newcommand{\abs}[1]{ \lvert #1 \rvert}

\newcommand{\td}[0]{\mathrm{d}}
\DeclareMathOperator{\support}{Supp}
\DeclareMathOperator{\var}{Var}
\begin{document}

\title{Reinforcement learning from comparisons: \\Three alternatives is enough, two is not.}
\author{Beno\^{\i}t \textsc{Laslier}\\Institut Camille Jordan, Universit\'e Claude Bernard Lyon 1,\\43 boulevard du 11 novembre 1918,\  69622 Villeurbanne, France
\and Jean-Fran\c{c}ois \textsc{Laslier}\\CNRS and \'{E}cole polytechnique,\\91128\ Palaiseau, France \\jean-francois.laslier@polytechnique.edu}
\maketitle

\begin{abstract}
The paper deals with the problem of finding the best alternatives on the basis of pairwise comparisons when these comparisons need not be transitive. In this setting, we study a reinforcement urn model.\ We prove convergence to the optimal solution when reinforcement of a winning alternative occurs each time after considering three random alternatives.
The simpler process, which reinforces the winner of a random pair does not always converges: it may cycle.

\end{abstract}


\section{Introduction}

In a pairwise comparison problem, one is given a set of alternatives, with data about how they compare the ones to the others.\ In its purest form, on which we focus in the present paper, we simply have, for any pair of distinct alternatives, the information of which one \textquotedblleft beats\textquotedblright\ the other.\ Such a data set is called a tournament.
Basic results on this structure can be found in Moon \cite{Moon68}.

For logical as well as practical reasons, binary relations are at the basis of choice theory. Two classical examples are sport competition and majoritarian politics.\ Many sports involve by definition two players (or teams), so that competition among any number of players must take the form of the aggregation of pairwise comparisons. In majority voting, a candidate is socially preferred to another if and only if a majority of the voters prefer the former to the latter.\ More generally, the prevalence of that kind of binary relations can be traced back to specific features of efficient natural languages (Rubinstein \cite{Rubinstein96}).

If a chess player beats all the other players, he or she is clearly the best.\ If a candidate cannot be defeated under majority rule by any challenger, that \textquotedblleft Condorcet\textquotedblright\ candidate can claim to be the best according to majority rule.\footnote{This observation does not imply that the majority principle is good for Politics.}\ But if no alternative beats all the others, it is not clear how to define the best alternatives.\ The problem of choosing from pairwise comparisons has thus attracted the attention of scholars in various fields, most often from the axiomatic, normative, point of view (David \cite{David63}; Fishburn \cite{Fishburn77}; Laslier \cite{Laslier97}%
; Brandt et al. \cite{Brandtetal12}).

In the present paper we tackle the same problem from an evolutionary perspective instead of an axiomatic one.\
We consider dynamic processes according to which, at each period in time, a small number (2 or 3) of alternatives are sampled, the tournament is played among these few alternatives, and the winning alternative is reinforced in the sense that it will be sampled with higher probability in the future.
Where does such a mechanical adaptive process go?
Using a standard urn model, where reinforcing an alternative is adding a colored ball to an urn, we obtain two results.\

(i) If one samples three alternatives (distinct or not) at each date, the process is able to discover the optimal solution of the tournament, that is the unique probability distribution $p^*$ which is, in expectation, defeated by no alternative.
With probability one, the composition of the sampling urn, which defines the probability $p_\tau$ of choosing the various alternatives at time $\tau$, tends to $p^*$ when $\tau$ tends to infinity.

(ii) If one samples only two alternatives at each date, the process is not able to discover the optimal solution, unless the solution is degenerated, with one alternative defeating all the others.
With probability one, the composition of the sampling urn, which gives the probability $p_\tau$ of choosing the alternatives, concentrates on the support of the optimal solution $p^*$.
This means that all the alternatives which are played with zero probability in the optimal solution are chosen with a probability going to $0$.
However the composition of the urn may cycle, staying away from the optimal solution. In some cases we even prove that it cycles with probability one.

The negative result (ii) echoes known results about the evolutionary instability of mixed equilibria in evolutionary game theory. For instance cycling with probability one is proven by Posch \cite{Posch97} in a reinforcement urn model for $2 \times 2$ games.\

The positive result (i) seems more original. In a study of imitation processes in Matching Pennies games, Hofbauer and Schlag \cite{HofbauerSchlag00} observe that players end up closer to the equilibrium if they sample several individual before imitating: there is still cycling, but closer to the equilibrium.\  Our results might be interpreted in the same spirit: learning slower leads to more stability.

The techniques we use to derive these results are standard in the field of adaptive processes with reinforcement; see Pemantle 2007 \cite{Pemantle07}, and they belong to the family of martingales techniques. The main ingredient of the proof is the definition of a well chosen function of the process whose values form a martingale (see \eqref{definition_mu}). We use the convergence theorem for positive martingales to obtain some global asymptotic information about the process. In the case of three alternatives we get fairly directly the convergence of the process while for the case of two alternatives the convergence theorem has to be complemented with a variance analysis to prove the non-convergence.

The paper is organized as follows. Section \ref{sec_framework} introduces the necessary notions about tournaments: definition and notation (\ref{subsec_tournaments}), the Markov chain induced by the play of small-size tournaments at each date (\ref{subsec_markov}), the tournament game which allows to define and to prove existence of the optimal solution (\ref{subsec_game}), and some further preliminary material (\ref{subsec_formulas}, \ref{subsec_relation}).\
Section \ref{sec_learning} starts by the definition of urns and of the adaptive processes (\ref{subsec_reinforcement}).
Then, in order to illustrate the argument in a simple way, a toy example is introduced and treated according to its deterministic approximation
(\ref{subsec_example}).
The statement and proof of our main result on three-alternatives reinforcement is found in (\ref{subsec_twostep}) and two-alternatives reinforcement is treated in (\ref{subsec_onestep}), before a short conclusion (\ref{subsec_conclusion}).

\section{Framework\label{sec_framework}}

\subsection{Tournaments\label{subsec_tournaments}}

Let $X$ be a finite set.\ A \emph{tournament} $T$ on $X$ is a complete and antisymmetric binary relation.\ For any $x$ and $y$ in $X$, one and only one of the three possibilities occurs: $x=y$, $x\,T\,y$, or $y\,T\,x$.\
When $x\,T\,y$ we often say that $x$ beats $y$.\ Define the sets :
\begin{equation}
T^{+}(x)=\{y\in X:xTy\},\ T^{-}(x)=\{y\in X:yTx\}. \label{def_ens_T}%
\end{equation}
The binary relation $T$ is fixed throughout this paper.\ It is sometimes
easier to use the notation:%
\[
\max\{x,y\}=\left\{
\begin{array}
[c]{c}%
x\\
y
\end{array}%
\begin{array}
[c]{l}%
\text{if}\ x\,T\,y\text{ or }x=y\\
\text{if }y\,T\,x\text{.}%
\end{array}
\right.
\]
An alternative which beats all other alternatives is called a \textit{Condorcet winner}. A tournament can have a Condorcet winner or not,
but cannot have two.\ The \emph{Top-Cycle} of the tournament $T$ is the smallest (by inclusion) set $Y\subseteq X$ such that:%
\[
\forall\ x\in X\setminus Y\text{, }\exists\ y\in Y:yTx.
\]
It is easily seen that such a set is unique and reduces to a singleton $\{c\}$ if and only if $c$ is a Condorcet winner. The literature on tournaments and formal political science has shown that the Top-Cycle is usually a very large set (McKelvey \cite{McKelvey79}), and has proposed many refinements of this set (see \cite{Laslier97} for a survey).

\subsection{A Markov chain\label{subsec_markov}}

Let $\Delta(X)$ be the set of probability distributions on $X$ and let $p\in\Delta(X)$. The support of $p$ is denoted by $\support(p)$. Given $p$, define a sequence $(p^{[t]})_{t\in%
\Nat
}$ of probability distributions on $X$ derived from $p$ in the following way :
\begin{align}
p^{[0]}  &  =p,\\
p^{[t+1]}(x)  &  =p^{[t]}(x)\cdot p\left(  T^{+}(x)\cup\{x\}\right)
+p^{[t]}(T^{+}(x))\cdot p(x), \label{pt+1}%
\end{align}
for any $t\geq0$, for any $x\in X$

The interpretation is that $p^{[t]}$ is the distribution of a random variable $\xi(t)\in X$ such that $\xi(0)$ is chosen at random according to $p$ and then, given that $\xi(t)=x$, $\xi(t+1)$ is the winner (according to $T$) of the comparison between $x$ and some alternative $y$ randomly chosen in $X$ according to $p$. Therefore $\xi(t+1)=x$ either because $\xi(t)$ was already equal to $x$ and $y$ was chosen in $T^{+}(x)\cup\{x\}$ (first term in the above formula), or because $\xi(t)$ was in $T^{+}(x)$ and $x$ was chosen according to $p$ (second term). We call $p$ the \textquotedblleft sampling\textquotedblright\ probability.\

This process is usually considered with $p$ uniform on $X$ (Daniels \cite{Daniels69}, Ushakov \cite{Ushakov76}, Levchenkov \cite{Levchenkov92}, Slutzky and Volij \cite{SlutzkiVolij06}, Chebotarev and Shamis \cite{ChebotarevShamis98, ChebotarevShamis06}).\
We need the general version because, later in this paper, $p$ will be endogenous.\

Given $p$, the stationary distribution for this finite Markov chain exists and is unique;\footnote{We state the results in this section without proofs.\ They are easily derived from elementary theory of finite Markov chains and have already been noticed for $p$ uniform in the mentioned references.} we denote it by $p^{[\infty]}$.\ It is characterized by the fact that $\support(p^{[\infty]})\subseteq\support(p)$ and, for any $x$ in
$\support(p)$,
\begin{equation}
p^{[\infty]}(T^{+}(x))\cdot p(x)=p^{[\infty]}(x)\cdot p(T^{-}(x)).
\label{pinfty}%
\end{equation}
Notice that the inclusion $\support(p^{[\infty]})\subseteq\support(p)$ may be strict; indeed,\ $p^{[\infty]}(x)=0$ when $p(T^{-}(x))=0$, that is when $x$ beats no alternative in the support of $p$. More exactly, $\support(p^{[\infty ]})$ is the Top-Cycle of the restriction of $T$ to $\support(p)$; thus $\support(p^{[\infty]})$ does not exactly really depends on $p$ but only on $\support(p)$.\ If $p$ has full support, for instance in the usual case where $p$ is uniform, $\support(p^{[\infty]})=TC(T)$.

\subsection{The tournament game\label{subsec_game}}

The tournament game is the two-player, symmetric, zero-sum game defined by the strategy set $X$ and the payoff function $g(x,y)=+1$ if $x\,T\,y$, $g(x,y)=0$ if $x=y$, and $g(x,y)=-1$ if $y\,T\,x$. For $p,q\in\Delta(X)$ two probability distributions on $X$, write:%
\begin{equation}
g(p,q)=\sum_{y\in X}g(x,y)p(x)q(y). \label{def_g}%
\end{equation}
From the definition, $g$ is clearly antisymmetric: $g(q,p)=-g(p,q)$.


The tournament game has been studied by graph theorists (Ficher and Ryan \cite{FisherRyan92, FisherRyan95a, FisherRyan95b}) and has more recently attracted attention of computer scientists (Rivest and Chen \cite{RivestChen10}). As a model of majority voting and two-party electoral competition, it studied in Social Choice theory and formal Political Science\ (Moulin \cite{Moulin86}, Myerson \cite{Myerson93, Myerson95}, Laslier \cite{Laslier00a, Laslier00b}).\ Remarkably, such a game has a unique equilibrium. Here is the precise result that will be needed in the sequel.
(Fisher and Ryan \cite{FisherRyan95a} prove this result using linear algebra and Laffond et al. \cite{LaffondetalGEB93} have a direct proof using a parity argument.)

\begin{proposition}
\label{definitionp*} There exists a unique $p^{\ast}\in\Delta(X)$ such that $g(p^{\ast},q)\geq0$ for all $q\in\Delta(X)$. This $p^{\ast}$, called the optimal strategy, is also characterized by the following : for all $x\in X$,
\begin{align*}
p^{\ast}(x)  &  >0\iff g(x,p^{\ast}(x))=0\\
p^{\ast}(x)  &  =0\iff g(x,p^{\ast}(x))<0.
\end{align*}
\
\end{proposition}

The support of the optimal strategy is called the \emph{Bipartisan Set} of the tournament: $\support(p^{\ast})=BP(T)$.\ This set is a subset of the Top Cycle and the inclusion is often strict.\ For instance, in totally random tournaments, the Top Cycle contains all the alternatives and the Bipartisan Set contains only half of them (Fisher and Reeves \cite{FisherReeves95}).

\subsection{Two formulas\label{subsec_formulas}}

Before we go further and explain the relation between the game optimal strategy and stationary probabilities, it is useful to state two technical formulas. The following lemma describes the probabilities $p^{[2]}$ and $p^{[3]}$, obtained after sampling two or three alternatives with the Markov chain defined in Section \ref{subsec_markov}, in term of the payoff function $g$.

\begin{lemma}
\label{lemmeformules} For any $x\in X$:$\;$%
\begin{align*}
p^{[2]}(x)  &  =p(x)\cdot(1+g(x,p)),\\
p^{[3]}(x)  &  =p(x)\cdot\left(  1+\tfrac{3}{2}g(x,p)+\tfrac{1}{2}%
g(x,p)^{2}+\sum_{y\in X}p(y)g(x,y)g(y,p)\right)  .
\end{align*}

\end{lemma}

\dem
First let us notice a useful equality. By definition \eqref{def_ens_T} and
\eqref{def_g}
\begin{equation}
g(x,p)=p(T^{+}(x))-p(T^{-}(x)),
\end{equation}
and, since $p(T^{+}(x))+p(T^{-}(x))+p(x)=1$, we get:
\begin{equation}
1+g(x,p)=2p(T^{+}(x))+p(x). \label{egalitedugain}%
\end{equation}

Let $a$ and $b$ be chosen according to $p$ and let $x=\max\{a,b\}$, then:
\begin{align*}
p^{[2]}(x)  &  =\Pr[a=x]\cdot\Pr[b\in T^{+}(x)\cup\{x\}]+\Pr[a\in
T^{+}(x)]\cdot\Pr[b=x]\\
&  =p(x)\cdot\left(  2p\left(  T^{+}(x)\right)  +p(x)\right) \\
&  =p(x)\cdot\left(  1+g(x,p)\right)  .
\end{align*}

For the second formula, let $a$, $b$ and $c$ be chosen according to $p$.\ An
alternative $x$ appears as $x=\max\left\{  \max\{a,b\},c\right\}  $ in the two exclusive cases:%
\begin{align*}
&  \left.  x=\max\{a,b\}\text{, and }x=\max\{x,c\}\right.  .\\
&  \left.  x\,T\,\max\{a,b\}\text{, and }x=c\right.  .
\end{align*}
In the first line, the event $x=\max\{a,b\}$ has probability $p^{[2]}%
(x)=p(x)\cdot(1+g(x,p))$ so the probability of the first case is $p(x)\cdot\left(  1+g(x,p)\right)  \cdot p(T^{+}(x)\cup\{x\})$.\ In the second
line, the event $x\,T\max\{a,b\}$ has probability
\[
p^{[2]}(T^{+}(x))=\hspace{-9pt}\sum_{y\in T^{+}(x)}\hspace{-9pt}%
p^{[2]}(y)=\hspace{-9pt}\sum_{y\in T^{+}(x)}\hspace{-9pt}p(y)(1+g(y,p)),
\]
therefore the probability $p^{[3]}$ is:%
\begin{align*}
p^{[3]}(x)  &  =p(x)\cdot\left(  1+g(x,p)\right)  \cdot\left(  p(x)+\hspace
{-9pt}\sum_{y\in T^{+}(x)}\hspace{-9pt}p(y)\right)  +p(x)\hspace{-9pt}%
\sum_{y\in T^{+}(x)}\hspace{-9pt}p(y)(1+g(y,p))\\
&  =p(x)^{2}\left[  1+g(x,p)\right]  +p(x)\hspace{-9pt}\sum_{y\in T^{+}%
(x)}\hspace{-9pt}p(y)\left[  \left(  2+g(x,p)\right)  +g(y,p)\right]  .
\end{align*}
Using the fact that $\frac{1+g(x,y)}{2}$ is $1$ if $y\in T^{+}(x)$, is $1/2$
if $y=x$, and is $0$ if not, one finds:%
\begin{align*}
&  \frac{p^{[3]}(x)}{p(x)}\\
&  =\sum_{y}p(y)\left[  \left(  2+g(x,p)\right)  +g(y,p)\right]
\frac{1+g(x,y)}{2}\\
&  =\tfrac{1}{2}\sum_{y}p(y)\left[  \left(  2+g(x,p)\right)
+g(y,p)+2g(x,y)+g(x,p)g(x,y)+g(y,p)g(x,y)\right] \\
&  =1+\tfrac{1}{2}g(x,p)+\tfrac{1}{2}g(p,p)+g(x,p)+\tfrac{1}{2}g(x,p)^{2}%
+\tfrac{1}{2}\sum_{y}p(y)g(y,p)g(x,y)\\
&  =1+\tfrac{3}{2}g(x,p)+\tfrac{1}{2}g(x,p)^{2}+\tfrac{1}{2}\sum
_{y}p(y)g(y,p)g(x,y),
\end{align*}
which is the announced formula. \QED

\subsection{Relation between optimal strategies and stationary probabilities\label{subsec_relation}}

We first observe that the game optimal strategy $p^{\ast}$ satisfies a nice
fixed-point property if we take $p^{[1]}=p^{\ast}$ as the sampling probability
to build the Markov chain, and that only an optimal strategy can be such a
fixed point.

\begin{proposition}
Let $p^{\ast}$ be the optimal strategy for the tournament game,
then\ $(p^{\ast})^{[2]}=(p^{\ast})^{[\infty]}=p^{\ast}$. Conversely, let $p$
be such that $p^{[2]}=p$, then $p$ is the optimal strategy for the tournament
game restricted to the support of $p$.
\end{proposition}

\dem
By lemma \ref{lemmeformules}, ${p^{\ast}}^{[2]}(x)=p^{\ast}(x)(1+g(x,p^{\ast
})$, and, by proposition \ref{definitionp*}, either $p^{\ast}(x)=0$ or $g(x,p^{\ast})=0$.

Conversely if $p^{[2]}(x)=p(x)=p(x)(1+g(x,p))$ then $g(x,p)=0$ as soon as $p(x)>0$ and $p$ is the optimal strategy on its support.
\QED

\section{Learning\label{sec_learning}}

With the previous background material in mind, we turn to the main result of this paper.\ Instead of considering re-sampling at each date according to a constant probability distribution, as is done in the previously described Markov chains, we describe learning processes where winning alternatives are reinforced at the level of the sampling probability.\ These processes can be implemented with random urns.

\subsection{Choice by reinforcement\label{subsec_reinforcement}}

An \emph{urn} on $X$ is a list $n$ of strictly positive integers $n(x)$, $x\in X$. The integer $n(x)$ is the \textquotedblleft number of balls of color $x$ in the urn $n$.\textquotedblright\ The set of such urns on $X$ is denoted by $\mathcal{N}$, formally:%
\[
\mathcal{N}=%
\Nat
_{+}^{X}.
\]
To each $n\in\mathcal{N}$ is associated the probability distribution
$\widetilde{n}$ on $X$ defined by%
\[
\widetilde{n}(x)=\frac{n(x)}{\sum_{y\in X}n(y)}.
\]
When we write that the alternative $x$ is picked in the urn $n$, we mean that
$x$ is picked in $X$ according to the probability $\widetilde{n}$.

A \emph{random urn sequence} is a sequence $U_{\tau}$, $\tau\in%
\Nat
$ of random variables on $\mathcal{N}$ such that $U_{\tau+1}$ is defined
conditionally on $U_{\tau}$. Here are three examples:

\begin{enumerate}
\item Two-alternatives reinforcement.\ Given a realization $n_{\tau}\in\mathcal{N}$ of
$U_{\tau}$, an alternative\ $x$ is picked in $X$ according to the probability
distribution $\widetilde{n_{\tau}}^{[2]}$, and one ball of color $x$ is added
to the urn: $n_{\tau+1}(w)=n_{\tau}(w)+1$ and for all $v\neq w$, $n_{\tau
+1}(v)=n_{\tau}(v)$.\ This means that two alternatives, say $a$ and $b$ are
picked independently in the urn $n_{\tau}$, and are compared according to $T$.
The result of the comparison is $x=\max\{a,b\}$, that is: $x=a$ if $a=b$ or if
$a\,T\,b$ and $x=b$ if $b\,T\,a$. Alternative $x$ is reinforced.

\item Three-alternatives reinforcement.\ Same thing as above, with the probability
distribution $\widetilde{n_{\tau}}^{[3]}$.\ This means that three
alternatives, say $a$, $b$ and $c$ are picked independently in $X$ according
to $n_{\tau}$; $a$, $b$ and $c$ are compared according to $T$ in sequence and
one ball of color $x=\max\left\{  \max\{a,b\},c\right\}  $ is added to the urn. Remark that there are only two cases : ranked alternatives where we reinforce the top one or a cycle where we reinforce at random.

\item Fast reinforcement. Same thing as above, with the probability
distribution $\widetilde{n_{\tau}}^{[\infty]}$, the stationary distribution
for $T$ when sampling is done according to $\widetilde{n_{\tau}}$.
\end{enumerate}

Remark that the first two examples can be concretely implemented easily, as
described, but fast reinforcement cannot.

\subsection{A motivational example\label{subsec_example}}

This section presents a simple non-rigorous argument to justify our focus on three-alternatives reinforcement.
Consider the simplest possible non trivial tournament : a cycle of three alternatives $A$, $B$ and $C$ with $A \,T\, B$, $B\,T\,C$, $C\,T\,A$. In order to evaluate the long term behavior of two-alternatives and three-alternatives reinforcement, we use a deterministic continuous time motion corresponding to the limit of a large number of balls in the urn. We write $a(t)$, $b(t)$ and $c(t)$ the ``number'' of balls of each type and $\tilde{a}(t)=a(t)/t$, $\tilde{b}(t)=b(t)/t$ and $\tilde{c}(t)=c(t)/t$ the corresponding probabilities.

For two-alternatives reinforcement we get :
\begin{eqnarray}
   a' &= &\tilde{a}^2 +2\tilde{a}\tilde{c} \nonumber \\
  b' & = & \tilde{b}^2 + 2\tilde{b}\tilde{a}  \\
  c' & = & \tilde{c}^2 + 2\tilde{c}\tilde{b} \nonumber
\end{eqnarray}
and we remark that
\begin{eqnarray*}
   \tfrac{ \td}{\td t}\left(\ln \tilde{a} + \ln \tilde{b} +\ln \tilde{c} \right) & = & \tfrac{ \td}{\td t}\left( -3\ln t + \ln a + \ln b + \ln c \right) \\
	& = &-\frac{3}{t}+ \frac{\tilde{a}^2 +2\tilde{a}\tilde{c}}{a}+ \ldots \\
	& = &-\frac{3}{t}+ \frac{\tilde{a}+2\tilde{c}}{t}+\frac{\tilde{b}+2\tilde{a}}{t}+\frac{\tilde{c}+2\tilde{b}}{t} \\
	& = & 0
\end{eqnarray*}
so $(\tilde{a},\tilde{b},\tilde{c})$ cannot converge to the optimal probability independently of the state at finite time.

It is then natural to study three-alternatives reinforcement, for which :
\begin{eqnarray}
   a' & = &\tilde{a}^3 +3\tilde{a}^2\tilde{c}+3\tilde{a}\tilde{c}^2+\tilde{a}\tilde{b}\tilde{c} \nonumber \\
  b' & = & \tilde{b}^3 + 3\tilde{b}^2\tilde{a}+3\tilde{b}\tilde{a}^2+ \tilde{a}\tilde{b}\tilde{c} \\
  c' & = & \tilde{c}^3 + 3\tilde{c}^2\tilde{b}+3\tilde{c}\tilde{b}^2+\tilde{a}\tilde{b}\tilde{c} \nonumber
\end{eqnarray}
and for the same quantity
\begin{eqnarray*}
   \tfrac{ \td}{\td t}\left(\ln \tilde{a} + \ln \tilde{b} +\ln \tilde{c} \right) & = & -\frac{3}{t} +\frac{ \tilde{a}^2+3\tilde{a}\tilde{c}+3\tilde{c}^2 +\tilde{b}\tilde{c}}{t} +\ldots \\
      & = & \frac{ 4(\tilde{a}^2+\tilde{b}^2+\tilde{c}^2)+4(\tilde{a}\tilde{c}+\tilde{b}\tilde{a}+\tilde{c}\tilde{b})-3}{t} \\
      & = & \frac{ \tilde{a}^2+\tilde{b}^2+\tilde{c}^2-2(\tilde{a}\tilde{c}+\tilde{b}\tilde{a}+\tilde{c}\tilde{b})}{t}
\end{eqnarray*}
and simple calculus shows that this last term is positive expect for $\tilde{a}=\tilde{b}=\tilde{c}=1/3$. Then $\ln \tilde{a} + \ln \tilde{b} +\ln \tilde{c}$ is an increasing negative function so it converges. It is not difficult to see, using the divergence of $\int 1/t \td t$, that this implies that $\tilde{a}^2+\tilde{b}^2+\tilde{c}^2-\tilde{a}\tilde{c}-\tilde{b}\tilde{a}-\tilde{c}\tilde{b}$ converges to $0$ and then that $(\tilde{a},\tilde{b},\tilde{c})$ converges to $(1/3,1/3,1/3)$ (the details of the arguments will be given in the rigorous proof of the next section).

With this example we can see that two-alternatives reinforcement should not converge to the optimal probability even for a simple tournament and when we neglect the effect of probabilistic noise while three-alternatives reinforcement seems to work in that case. In the next section we will prove that three-alternatives reinforcement actually converges for any tournament. We will use the same idea of computing the variation of $\ln \tilde{a} + \ln \tilde{b} +\ln \tilde{c}$ with technical changes for the general tournament, the discrete time and the probabilistic evolution.

\subsection{Three-alternatives reinforcement and martingale technique\label{subsec_twostep}}

We will now prove the result about three-alternatives reinforcement:

\begin{theorem}\label{thm:3alternatives}
For any initial urn $n_{0}\in\mathcal{N}$, the random urn sequence obtained by
three-alternatives reinforcement is such that the realization $n_{\tau}$, $\tau\in%
\Nat
$ almost surely verifies:
\[
\lim_{\tau\rightarrow\infty}\widetilde{n_{\tau}}=p^{\ast}\text{.}%
\]
\end{theorem}

The same proof will actually also give the first part of result about two-alternatives reinforcement, which we thus state now:
\begin{theorem}\label{thm:2alternatives_v1}
 For any initial urn $n_{0}\in\mathcal{N}$, the random urn sequence obtained by
two-alternatives reinforcement is such that the realization $n_{\tau}$, $\tau\in%
\Nat
$ almost surely verifies:
\[
\forall x\in X, p^\ast(x) = 0 \Rightarrow \lim_{\tau\rightarrow\infty}\widetilde{n_{\tau}}(x)=0\text{.}%
\]
\end{theorem}

The proof relies mainly on the study of a well chosen function of the state of the urn. Let $LD$ denotes a discrete logarithm: for integers $0<a<b$,%
\begin{equation}
LD[a,b]=-\sum_{i=a}^{b-1}\frac{1}{i}. \label{definition_LD}%
\end{equation}
Recall that at time $\tau\in%
\Nat
$, $n_{\tau}(w)$ denotes the number of $w$-balls in the urn. The total number of balls is increasing by $1$ at each time, so $\sum_{w}n_{\tau}(w)=A+\tau$. The probability of drawing a $w$-ball is $\tilde n_{\tau}(w)=n_{\tau}(w)/(A+\tau)$.
Consider the quantity%
\begin{equation}
\mu_{\tau}=\sum_{w\in X}LD[n_{\tau}(w),A+\tau]\cdot p^{\ast}(w), \label{definition_mu}
\end{equation}
that is the expected value, according to the optimal probability $p^{\ast}$, of the discrete logarithm at time $\tau.$

\begin{proposition}
\label{mu_est_ssmartingale} For both two-alternatives and three-alternatives reinforcement, the sequence $\mu_{\tau},$ $\tau\in%
\Nat
$ is a negative sub-martingale. More precisely we have, for two-alternatives reinforcement:
\[
\mathrm{E}\left[  \mu_{\tau+1}-\mu_{\tau}\mid n_{\tau}\right] = \frac{g(p^\ast,\tilde n_\tau )}{A+\tau}.
\]
and for three-alternatives reinforcement:
\begin{align*}
& \mathrm{E}\left[  \mu_{\tau+1}-\mu_{\tau}\mid n_{\tau}\right] \\
& = \frac{1}{A+\tau} \left(\tfrac{1}{2}g(p^{\ast},\tilde n)+\tfrac{1}{2}\sum_{w\in X}g(w,\tilde n)^{2}p^{\ast
}(w)+\sum_{v}\tilde n(v)g(p^{\ast},v)(1+g(v,\tilde n)) \right). \\
\end{align*}
\end{proposition}

\dem
We will write $p$ for $\tilde n_\tau$ and let $i$ denote either $2$ or $3$.
From $\tau$ to $\tau+1$, one and only one ball is added.\ This ball has type
$w$ with probability $p^{[i]}(w)$. Thus:
\begin{align*}
&  \mathrm{E}\left[  \mu_{\tau+1}-\mu_{\tau}\mid n_{\tau}\right]  \\
&  =-\sum_{w\in X}p^{[i]}(w)\left(  \sum_{v\neq w}\frac{1}{A+\tau}\cdot
p^{\ast}(v)+\left[  \frac{1}{A+\tau}-\frac{1}{n_{\tau}(w)}\right]  \cdot
p^{\ast}(w)\right)  \\
&  =\frac{-1}{A+\tau}+\sum_{w\in X}p^{[i]}(w)\frac{1}{n_{\tau}(w)}p^{\ast}(w)\\
& = \frac{-1}{A+\tau}+\sum_{w\in X}\frac{p^{[i]}(w)}{p(w)}\frac{p^{\ast}(w)}{A+\tau}\\
\end{align*}
Where in the last line we used the definition $p(w)=n_{\tau}(w)/(A+\tau)$. Using the formula for $p^{[i]}$ of lemma \ref{lemmeformules}, it comes, for two-alternatives:%
\begin{align*}
&  \mathrm{E}\left[  \mu_{\tau+1}-\mu_{\tau}\mid n_{\tau}\right]  \\
& = \frac{-1}{A+\tau}+\sum_{w\in X} \left( 1 + g(w,p) \right) \frac{p^{\ast}(w)}{A+\tau} \\
& = \frac{ g(p^\ast,p)}{A+\tau}.
\end{align*}
which is always non-negative. Furthermore, $g(p^{\ast},p)=0$ implies, by the first part of proposition~\ref{definitionp*}, that $\support(p)\subseteq\support(p^{\ast})$.

For three-alternatives, we have:
\begin{align*}
&  \mathrm{E}\left[  \mu_{\tau+1}-\mu_{\tau}\mid n_{\tau}\right]  \\
&  =\frac{-1}{A+\tau}+\sum_{w\in X}\left(  1+\tfrac{3}{2}g(w,p)+\tfrac{1}%
{2}g(w,p)^{2}+\sum_{v}p(v)g(w,v)g(v,p)\right)  \frac{p^{\ast}(w)}{A+\tau}.
\end{align*}
one can re-arrange:%
\begin{align*}
&  (A+\tau)\mathrm{E}\left[  \mu_{\tau+1}-\mu_{\tau}\mid n_{\tau}\right]  \\
&  =\tfrac{3}{2}g(p^{\ast},p)+\tfrac{1}{2}\sum_{w\in X}g(w,p)^{2}p^{\ast
}(w)+\sum_{v}p(v)g(p^{\ast},v)g(v,p)\\
&  =\tfrac{1}{2}g(p^{\ast},p)+\tfrac{1}{2}\sum_{w\in X}g(w,p)^{2}p^{\ast
}(w)+\sum_{v}p(v)g(p^{\ast},v)(1+g(v,p)).
\end{align*}
And all the terms in this sum are non-negative.\ The sum can be $0$ only if
both $\support(p)\subseteq\support(p^{\ast})$,
and $g(w,p)=0$ for all $w$ in the support of $p^{\ast}$, which implies
$p=p^{\ast}$ by the uniqueness in Proposition~\ref{definitionp*}.
\QED

\bigskip

We are now able to prove the two results of the beginning of the section.

\noindent\textbf{Proof of Theorems \ref{thm:3alternatives} and \ref{thm:2alternatives_v1}.}
We consider for this proof either two or three alternatives reinforcement. We have:
\begin{align}
\mathbb{E}[\mu_{\tau}]  &  =\mathbb{E}[\mu_{0}]+\mathbb{E}\left[  \sum
_{t=1}^{\tau}\mu_{t}-\mu_{t-1}\right] \\
&  =\mathbb{E}[\mu_{0}]+\mathbb{E}\left[  \sum_{t=1}^{\tau}\mathbb{E}\left[
\mu_{t}-\mu_{t-1}|\mu_{t-1}\right]  \right]  .
\end{align}
By Proposition \ref{mu_est_ssmartingale}, $\mu_{\tau}$ is a negative sub-martingale, so it converges almost surely
to an integrable random variable $\mu_{\infty}$ (see Corollary VII.4.1 and VII.4.2 in \cite{Shiryaev}).
Furthermore, the right hand side is an increasing negative sequence so it converges to a finite value. In the left hand side, the sum is an increasing
sequence of positive random variables so by the monotonous convergence theorem (Theorem II.6.1 in \cite{Shiryaev}) we can take the limit inside the expectation. Hence $\mathbb{E}\left[  \sum_{t=1}^{\infty}\mathbb{E}[\mu_{t}-\mu_{t-1}|\mu_{t-1}]\right]  =\lim\mathbb{E}[\mu_{\tau}]-\mathbb{E}%
[\mu_{0}]$ is finite and so $\sum_{t=1}^{\infty}\mathbb{E}[\mu_{t}-\mu_{t-1}|\mu_{t-1}]$ is almost surely finite.


Let $f^{[2]}(p) = g(p^\ast,p)$ and $f^{[3]}(p)=g(p,p^{\ast})+\sum_{w\in X}g(w,p)^{2}p^{\ast}(x)$.
The simplex $\Delta(X)$ is embedded in ${\mathrm{I\hspace{-0.15em}R}}^{X}$ so we use the
$L_{\infty}$ distance. With this distance, $f^{[i]}$ is continuous and
$d(\widetilde{n_{\tau}},\widetilde{n_{\tau+1}})\leq\frac{1}{A+\tau}$ almost surely.
Denote by $B(p,\eta)$ the ball of center $p$ and radius $\eta$.

Now consider a single realization of the urn process. Since $X$ is a finite
set, $\Delta(X)$ is compact, so let $\widetilde{n_{\infty}}$ be an accumulation
point for $\widetilde{n_{\tau}}$. We will show that necessarily
$f^{[i]}(\widetilde{n_{\infty}}) =0$. Looking for a contradiction, suppose
$f^{[i]}(\widetilde{n_{\infty}}) >0$.

Since $f^{[i]}$ is continuous, let $\epsilon,\eta>0$ be such that $\forall p\in B(\widetilde{n_{\infty}},\eta)$,
$f^{[i]}(p)>\epsilon$. Let $\phi$ be a sub-sequence such that $\forall
\tau,\ \widetilde{n_{\phi(\tau)}}\in B(\widetilde{n_{\infty}},\eta/2)$ and
$\phi(\tau+1)>(1+\eta)\phi(\tau)$. Then:
\begin{align}
\sum_{t=0}^{\infty}\mathbb{E}[\mu_{t+1}-\mu_{t}|n_{t}]  &  \geq\sum_{t=0}^{\infty}\frac{1}{2(A+t)}f^{[i]}(\widetilde{n_{t}})\\
&  \geq\sum_{\tau=0}^{\infty}\frac{1}{2(A+\phi(\tau))}\sum_{t=\phi(\tau
)}^{\lfloor(1+\eta/2)\phi(\tau)\rfloor}f^{[i]}(\widetilde{n_{t}})\\
&  \geq\sum_{\tau=0}^{\infty}\frac{1}{2(A+\phi(\tau))}\lfloor(1+\eta
/2)\phi(\tau)\rfloor\epsilon
\end{align}
The right hand side of the last line is infinite. We already proved that the left hand side is almost surely finite.\ This contradiction proves that all the accumulation points of $\widetilde{n_{t}}$ are zeros of $f^{[i]}$.\ It follows that $f^{[i]}(\widetilde{n_{t}})\rightarrow 0$ almost surely. We have seen in the proof of Proposition \ref{mu_est_ssmartingale} that this fact implies exactly the theorems.
\QED

\subsection{Two-alternatives reinforcement and variance estimates\label{subsec_onestep}}

In this section, we study in detail the two-alternatives reinforcement. The main idea is to study the variance of $\mu_\infty$ conditionally on the state at a large time $t$.


\begin{theorem}\label{th:two_step}
For any initial urn $n_0 \in \mathcal{N}$, the random urn sequence obtained by two-alternatives reinforcement is such that:
\begin{enumerate}
 \item \label{th:two_step:pt1} almost surely, for all alternatives $x$ such that $p^\ast(x)=0$, $\tilde{n}_\tau (x) \rightarrow 0$ when $\tau \rightarrow \infty$;
\item \label{th:two_step:pt2} with positive probability, the realized sequence $\tilde{n}_\tau, \tau \in \Nat$ has no limit as $\tau \rightarrow \infty$;
\item if $T$ is such that $\forall x, p^\ast(x) > 0$ (in other words, $BP(T) = X$) then, with probability one, $\tilde{n}_\tau, \tau \in \Nat$ has no limit.
\end{enumerate}
\end{theorem}

To simplify notation, in this section we let the process start at $\tau \neq 0$ so that $\tau$ always denote the number of ball in the urn (i.e. $A=0$). We will also only consider two-alternatives reinforcement in this section. Recall the piece of notation $\support(p^\ast) = BP$. Also recall from the last section that $\mu_\tau$ is a negative submartingale so it has an almost sure limit $\mu_\infty$. Let $\phi = \sum_{x\in BP} p^\ast(x) \log p^\ast(x)$ be the value of $\mu_\infty$ when $\tilde n_{\tau}$ converges to $p^\ast$.

The first point is the following variance estimate :
\begin{lemma}\label{lemme:variance}
Let $\tau > 0$ and let $\epsilon_\tau(x) = p^*(x)/\tilde{n}_\tau(x)-1$. We have
\begin{equation}
\mathbb{E}[ (\mu_{\tau+1}-\mu_{\tau})^2 | \mathcal{F}_\tau ] = \frac{1}{\tau^2} \sum_{x\in BP} \tilde{n}^{[2]}(x) \epsilon_\tau(x)^2
\end{equation}
\end{lemma}
\dem
This is a straightforward computation :
\begin{align}
    \mathbb{E}[ (\mu_{\tau+1}-\mu_{\tau})^2 | \mathcal{F}_\tau ]  & = \sum_x \tilde n_\tau^{[2]} (x) \bigl(-\sum_y p^*(y)\frac{1}{\tau} + p^*(x)\frac{1}{n(x)} \bigr)^2 \\
                & = \sum_x \tilde n_\tau^{[2]} (x) \frac{1}{\tau^2}( -1 + \frac{p^*(x)}{\tilde n(x)})^2 \\
                & = \frac{1}{\tau^2} \sum_x \tilde n^{[2]}(x) \epsilon_\tau(x)^2
\end{align}
\QED

The main point is the factor $\frac{1}{\tau^2}$ that make the series of those terms summable (once a small control on $\epsilon$ is provided). Thus the variance of $\mu_\infty$ conditioned on $\mathcal{F}_\tau$ will be of order $\epsilon^2/\tau$. Thus with hight probability $\mu_\infty$ will be close to $\mu_\tau$ so that if $\mu_\tau$ is far enough from $\phi$, then $\mu_\infty \neq \phi$.

We will first consider the case where $BP \neq X$.
In this case we have $\mathbb{E} [ \mu_{\tau+1}-\mu_\tau \mid \mathcal F_\tau] = \frac{g (p^\ast, \tilde n_{\tau})}{\tau} \geq g_0 \cdot \tilde{n} (BP^c) >0 $ (where $g_0 = \inf_{y\in BP^c} g(p^\ast,y)$) so we need an estimate of $\tilde n(BP)$.

\begin{lemma}\label{lemme:control_drift}
 Suppose that there exists $\pi \in (0,1)$ such that, at each time $\tau$ the probability of adding a ball in $BP^c$ is inferior to $\pi \cdot \tilde{n}_\tau (BP^c)$. Then we have for all $\tau \geq \tau_0$,
\begin{align*}
 \mathbb E [ \tilde{n}_\tau(BP^c) \mid \mathcal F_\tau ] &\leq \prod_{t=\tau_0+1}^{\tau} \frac{t+\pi-1}{t} \mathbb E[\tilde n_{\tau_0} ] \\
    & \leq \left( \frac{\tau_0}{\tau} \right)^{1-\pi} \tilde n_{\tau_0}.
\end{align*}
\end{lemma}

\dem
The first line comes from a straightforward induction
\begin{align*}
  \mathbb{E} [\tilde n_\tau ] &= \tfrac{1}{\tau}\mathbb E[ n_{\tau-1}] + \tfrac{1}{\tau} \mathbb E[ n_\tau -  n_{\tau-1}]\\
      & =\tfrac{\tau-1}{\tau} \mathbb E[ \tilde n_{\tau-1} ] + \tfrac{1}{\tau} \mathbb E\bigl[ \mathbb E[n_\tau -  n_{\tau-1} \mid \mathcal F_{\tau-1}] \bigr] \\
      & \leq \tfrac{\tau-1}{\tau} \mathbb E[ \tilde n_{\tau-1} ] + \tfrac{1}{\tau} \mathbb E[ \pi\tilde n_{\tau-1} ] \\
      & \leq \tfrac{\tau+\pi-1}{\tau} \mathbb E[ \tilde n_{\tau-1} ]
\end{align*}
and for the second line we have
\begin{align*}
 \log \prod_{t=\tau_0+1}^{\tau} \frac{t+\pi-1}{t} & \leq \sum_{\tau_0+1}^\tau \log( 1+\frac{\pi-1}{t} )\\
    & \leq \sum_{\tau_0+1}^\tau \frac{\pi-1}{t} \\
    & \leq (\pi-1)\log(\frac{\tau}{\tau_0})
\end{align*}
\QED

Now we are able to prove the third point of the theorem.

\textbf{Proof of theorem \ref{th:two_step}, case $BP^c \neq \emptyset$.}
 Remark that for $\delta >0$ small enough and $\tau_0$ big enough, the set
 \begin{equation*}
 S=\{ n \in \mathcal N \text{ s.t } n(X) \geq \tau_0 \text{ and } \abs{\phi - \sum_x p^\ast(x) LD(n(x), n(X))} \leq \delta \}
 \end{equation*}
 verifies :
\begin{itemize}
 \item $\forall n\in S, \forall x \in BP, \abs{ p^\ast(x)/\tilde n(x) -1} \leq 1$
  \item $\forall n \in S, \forall x \in BP^c, 1+g(x,\tilde n) \leq \pi$
\end{itemize}
for some $\pi > 1 - \inf_{x\in BP^c} g(p^\ast,x)$.

Let us assume that $\tilde{n}_{\tau_0} \in S$, (which clearly happens with positive probability). Let $T$ be the first time after $\tau_0$ such that $\tilde n_\tau \notin S$. $T$ is a stopping time so $\mu_{\tau \wedge T}$ is still a submartingale, let us call it $\mu'_\tau$. Furthermore, up to time $T$, we have
\begin{align*}
\mathbb E [\mu_{\tau+1} - \mu_\tau \mid \mathcal F_\tau] &= \frac{g(p^*,\tilde n_\tau)}{\tau}\\
    &\leq \frac{\tilde n_\tau(BP)}{\tau},\\
\end{align*}
and thus
\begin{align*}
 \mathbb E [\mu'_{\tau + 1} - \mu'_{\tau} \mid \mathcal F_{\tau_0}] & \leq \frac{1}{\tau} \mathbb E[ \tilde n_\tau(BP) 1_{\{ T \leq \tau \}} \mid \mathcal F_{\tau_0} ] \\
 & \leq \tau_0^{1-\pi} \tau^{pi-2} \tilde n_{\tau_0} (BP).\\
\end{align*}

Now let $\mu'_\infty$ denote the almost sure limit of $\mu_{\tau \wedge T}=\mu'_\tau$. Remark that if $\abs{\phi - \mu_\infty} (\omega) < \delta$ then $T(\omega) = \infty$ and thus $\mu_\infty (\omega) = \mu'_\infty(\omega)$. It is therefore enough to show that, with positive probability $\phi-\delta < \mu'_\infty < \phi$.

Note that $\mu'$ is a bounded submartingale, so it also converges in $L^1$ and $L^2$ toward $\mu'_\infty$. Thus

\begin{align*}
 \mathbb E [\mu'_\infty - \mu'_{\tau_0} \mid \mathcal F_\tau ] & = \sum_{t=\tau_0}^\infty \mathbb E[ \mu'_{t+1} - \mu'_t \mid \mathcal F_{\tau_0}] \\
      & \leq \sum_{t=\tau_0}^ \infty \tau_0^{1-\pi} t^{pi-2} \tilde n_{\tau_0} (BP) \\
      & \leq C \tilde n_{\tau_0}(BP),
\end{align*}
and
\begin{align*}
 \var (\mu'_\infty \mid \mathcal F_{\tau_0}) &\leq \sum_{t=\tau_0}^\infty \var( \mu'_{t+1} \mid \mathcal F_{\tau_0} ) \\
      & \leq \sum_{t=\tau_0}^\infty \mathbb E[ (\mu'_{t+1} - \mu'_t)^2 \mid \mathcal F_{\tau_0} ] \\
      & \leq \sum_{t=\tau_0}^\infty \frac{1}{t^2} \sum_{x\in BP} \tilde n^{[1]}(x) \epsilon_t (x)^2 \\
      & \leq C \frac{1}{\tau_0}. \\
\end{align*}

Finally consider a $\tau_0$ large enough so that $\frac{1}{\tau_0} \ll \delta$. It is clear that with a positive probability, $\tilde n_{\tau_0} \ll \delta$ and $\abs{\mu_{\tau_0} -\phi}$ is close to $\delta/2$. Under this event we see that $\abs{\phi-\mu'_\infty}$ is a random variable with expectation close to $\delta/2$ and variance small with respect to $\delta$ so $\abs{\mu'_\infty - \phi}$ has a positive probability to be in $[\delta/4,3\delta/4]$ which proves the theorem.
\QED

\bigskip
Now we turn to the case where $X = BP$. The idea will be similar, with Lemma \ref{lemme:variance} being the core argument. The main simplification comes from the fact that in this case $g(p^*,p) = 0$ for all $p$ so $\mu$ is a martingale and Lemma \ref{lemme:control_drift} will no longer be needed. However, in order to prove that $\mu_\infty$ is almost surely different from $\phi$, we will need an almost sure lower bound on $\abs{\phi-\mu_\tau}$ which will come from a careful analysis of the difference between discrete and real logarithms.
Finally since the almost sure bound that we will get will be much worse than the one we were able to have with positive probability, we will need to be more careful in our use of Lemma \ref{lemme:variance}.

First recall the following well known approximation result :
\begin{proposition}
Let $k>0$, there exists a constant $\gamma$ (Euler's constant) such that
\begin{equation}
    \log(k+1)+\gamma -\frac{1}{2}\sum_{i=k+1}^{\infty} \frac{1}{i^2} \leq  \sum_{i=1}^k \frac{1}{i} \leq \log(k+1)+\gamma -\frac{1}{2}\sum_{i=k+2}^{\infty} \frac{1}{i^2}
\end{equation}
Furthermore when $k$ tends to infinity,
\begin{equation}
\sum_{i=k+1}^{\infty} \frac{1}{i^2} \sim \frac{1}{k}
\end{equation}
\end{proposition}

This proposition implies the following corollary:
\begin{corollary}\label{lemme:distance_deterministe}
Let $T$ be a tournament on the set $X$ such that $BP=T$. There exists $c>0$ such that, for any urn $n \in \mathcal N$, 
\begin{equation}
   \sum_x p^*(x) LD( n(x), \tau ) \leq \phi - \frac{c}{\tau}.
\end{equation}
Furthermore, writing $\epsilon(x) = p^*(x)/\tilde{n}(x) -1$, if we restrict ourselves to large enough $\tau$ (with $\epsilon$ staying bounded) the constant $c$ can be taken as close as we want to :
\begin{equation}
 \frac{ \abs{X} -1 + \sum_{x\in BP} \epsilon(x)}{2}
\end{equation}
\end{corollary}
\dem
This is a straightforward computation using the definition of $LD$.
\QED

\bigskip
We also need a control of $\epsilon$ in term of $\mu$
\begin{lemma}\label{lemme:controle_epsilon}
Almost surely, for any time $\tau$
\begin{equation}
\sum_{x\in BP} [\tilde n_\tau(x) + p^*(x)/2]\epsilon_\tau(x)^2 \leq \phi - \mu_\tau
\end{equation}
\end{lemma}
\dem
We have
\begin{align}
\phi - \mu_\tau & \geq \sum_{x \in BP} p^*(x)[ \log p^* - \log p ] \\
                & \geq \sum_{x \in BP} p^*(x) ( \epsilon_\tau(x) + \epsilon_\tau(x)^2/2) \\
                & \geq \sum_{x \in BP} p^*(x) \epsilon_\tau(x)^2/2 + \sum_{x \in BP} \tilde n_\tau(x)(1+\epsilon_\tau(x))\epsilon_\tau(x) \\
                & \geq \sum_{x \in BP} (p^*(x)/2 +  \tilde n_\tau(x)) \epsilon_\tau(x)^2
\end{align}
because $\sum_{x \in BP} \tilde n_\tau(x) \epsilon(x) = \sum_{x \in BP} p^*(x) -\tilde n_\tau(x)=0$.
\QED

\bigskip
Together, the last lemmas have the following consequence
\begin{lemma}\label{lemme:reste_loin}
Let $\tau_0 > 0$ large enough such that $\abs{\phi-\mu_{\tau_0}} \geq \frac{5}{\tau_0-1}$, then there exists $\pi > 0$ such that, with probability at least $\pi$,
\[
        \forall \tau > \tau_0, \abs{\phi-\mu_{\tau}} \geq \frac{1}{\tau_0}.
\]
\end{lemma}
\dem
Let $d= \abs{\phi- \mu_{\tau_0}}$ and let $T= \inf\{ \tau > \tau_0 | \phi-\mu_{\tau} \geq 2d \text{ or } \phi-\mu_{\tau} \leq d/5 \}$. $T$ is a stopping time so $\mu'_\tau=\mu_{\tau \wedge T}$ is still a sub-martingale ; by definition it is also bounded so it converges almost surely and in all $L^p$. Let $\mu'_\infty$ denote its limit.
We will show that $\mathbb{P}( \mu'_\infty \in (d/5,2d) )> \pi >0$ for some $\pi$.
Remark that, since $\mu$ makes vanishing steps, we are only interested in the behaviour of $\mu$ close to $\phi$, so we can restrict ourself to small $\epsilon$.

As long as $\tau < T$, by definition we have $\phi-\mu_{\tau} \leq 2d$ and thus by Lemma \ref{lemme:controle_epsilon}, $\sum_x (p^*(x)/2 +  \tilde n_\tau(x)) \epsilon_\tau(x)^2 \leq 2d$. For $\epsilon$ small enough, this implies $\sum_x \tilde n_\tau^{[1]} \epsilon(x)^2 \leq 2d$ and thus, by Lemma \ref{lemme:variance}:
\begin{equation}
    \forall \tau > \tau_0, \mathbb{E} [ (\mu'_{t+1} - \mu'_{t})^2 | \mathcal{F} ] \leq \frac{2d}{\tau^2}.
\end{equation}
Summing up to infinity (recall that $\mu'$ has constant expectation):
\begin{align}
    \var(\mu'_\infty - \mu'_{\tau_0} | \mathcal{F}_{\tau_0} ) & = \sum_{\tau = \tau_0}^{\infty} \mathbb{E} [ ( \mu'_{\tau+1}-\mu'_\tau)^2 | \mathcal{F}_{\tau_0} ] \\
      & \leq \sum_{\tau = \tau_0}^{\infty} \frac{2d}{\tau^2} \\
      & \leq 2d \frac{1}{\tau_0-1} \\
      & \leq \frac{2}{5} d^2
\end{align}
where we used the hypothesis $d \geq \frac{5}{\tau_0-1}$ in the last line.

Remark that, since $\mu_{T \wedge \tau}$ is bounded, $\mathbb E [\mu'_\infty] = d$.
Moreover notice that a random variable with expectation $d$ which never takes value in the interval $(d/2,2d)$ has at least variance $d^2/2$.
Since  $\var (\mu'_\infty | \mathcal F_\tau) \leq \frac{2}{5} d^2$, we have $\mathbb{P} (\mu' \in (d/2,2d) ) \geq 1/10$ and on this event, $T=\infty$ so $\abs{\phi-\mu}$ has never reached $d/2 \geq \frac{1}{\tau_0}$.
\QED
\bigskip

\noindent\textbf{Proof of Theorem \ref{th:two_step}, case $BP = \emptyset$.}
First note that if $\tilde n_\tau$ converges, it has to be toward a fixed point.
It is easy to see that the fixed points of the dynamics are exactly the optimal strategies $p^*_Y$ corresponding to all subtournaments $Y \subseteq X$. (This of course includes $p^*=p^*_X$ itself.)
For $Y \subsetneq X$, $\sum_x p^*_X(x) \log p^*_Y(x) = -\infty$ so by the using the Markov inequality on $\mu$ we see that convergence to those fixed points is impossible.
Therefore we only have to rule out convergence towards $p^*_X$.

We first consider the case $\abs{X} \geq 12$. Then by Corollary \ref{lemme:distance_deterministe}, for any $\tau$ large enough, $\mu_\tau$ almost surely verifies the hypothesis of Lemma \ref{lemme:reste_loin}.
Fix any suitable $\tau_0$, by Lemma \ref{lemme:reste_loin}, $\mathbb{P} (\tilde n_\tau \rightarrow p^*) \leq \mathbb{P}( \exists \tau \geq \tau_0 | \abs{\phi-\mu_\tau} \leq 1/\tau_0) \leq 1-\pi$.
On the event that $\abs{\phi-\mu_\tau}$ does reach $1/\tau_0$ at time $\tau_1$, we can use Lemma \ref{lemme:reste_loin} at time $\tau_1$ to get $\mathbb{P} (\tilde n_\tau \rightarrow p^*) \leq (1-\pi)^2$.
By induction we get $\mathbb{P} (\tilde n_\tau \rightarrow p^*) = 0$ which proves the theorem.

For the case $3 \leq \abs{X} \leq 11$ (a non trivial tournament has at least $3$ elements), consider any $\tau_0$ large enough. By Corollary \ref{lemme:distance_deterministe} we have $\abs{\phi -\mu_{\tau_0}} \geq 1/\tau_0$. Let $T = \inf \{ \tau >\tau_0 | \abs{\phi -\mu_{\tau}} \leq 1/2\tau_0 \text{or} \abs{\phi -\mu_{\tau}} \geq \frac{5}{\tau_0}$. The event $\{ T= \infty \text{ or } \abs{\phi -\mu_{T}} \geq \frac{5}{\tau_0} \}$ has probability at least $1/9$ and if $\abs{\phi -\mu_{T}} \geq \frac{5}{\tau_0}$ we can apply Lemma \ref{lemme:reste_loin} at time $T$ so the conclusion of Lemma \ref{lemme:reste_loin} is still true with $1/\tau_0$ replaced by $1/2\tau_0$ and we can use the same induction as before to prove the theorem.
\QED

\subsection{Conclusion\label{subsec_conclusion}}

We have found the behavior of learning process designed to discover the \textquotedblleft best\textquotedblright\ alternatives in a tournament. Learning is achieved through the following idea. An alternative which is considered as \textquotedblleft good\textquotedblright\ at some date is reinforced for the future in the sense that one (slightly, and less and less) increases the probability for this alternative to be considered: reinforcement updates the sampling, or \textquotedblleft prior\textquotedblright\ probability. The test according to which an alternative is considered as a good one at time $t$ rests on comparing a few randomly chosen alternatives.

We found a very different behavior between the processes where reinforcement occurs after sampling two or three alternatives.
With three alternatives, the process converges almost surely to a well-defined limit that has a nice interpretation in term of the tournament game: it is the optimal strategy for this game.
One can therefore say that this form of learning is \textquotedblleft successful \textquotedblright.
With two alternatives, the picture is more complicated. The learning process \textquotedblleft succeeds \textquotedblright\ in finding the Bipartisan set (a set which has been argued to be more important in term of social choice than the numerical values of the optimal probabilities \cite{Laslier00a}), but not the optimal probabilities themselves.
We conjecture that the almost sure non-convergence happens for all tournaments and not only when $BP=X$.

\section*{Acknowledgment}

Thanks to Bastien Mallein for useful discussions in the early stage of this study.


\begin{thebibliography}{99}                                                                                               %


\bibitem {Brandtetal12}Brandt,Felix, Maria Chudnovsky, Ilhee Kim, Gaku Liu,
Sergey Norin, Alex Scott, Paul Seymour, and Stephan Thomasse (2011) A
counter-example to a conjecture of Schwartz.\ \textit{Social Choice and
Welfare} forthcoming.

\bibitem {ChebotarevShamis98}Chebotarev, P.\ T., Shamis, E. 1998.
Characterizations of scoring methods for preference aggregation.
\textit{Annals of Operation Research} \textbf{80}: 299---332.

\bibitem {ChebotarevShamis06}Chebotarev, P.\ T., Shamis, E. 2006. Preference
fusion when the number of alternatives exeeds two: indirect scoring
procedures. arXiv:math/060217v3 [math.OC]

\bibitem {Daniels69}Daniels, H. E. 1969. Round-robin tournament scores.
\textit{Biometrika} \textbf{56: }295---299.

\bibitem {David63}David, H. 1963. \textit{The Method of Paired Comparisons},
Griffin's statistical monographs and courses. Charles Griffins, London.

\bibitem {Fishburn77}Fishburn, P.C. 1977. \textquotedblleft Condorcet social
choice functions\textquotedblleft\ \textit{SIAM Journal on Applied
Mathematics} 33: 469---489.

\bibitem {FisherReeves95}David C. Fisher and Richard B. Reeves, 1995. Optimal
strategies for random tournament games.\ \textit{Linear Algebra and its
Applications}, \textbf{217}: 83---85.

\bibitem {FisherRyan92}Fisher, D., Ryan, J. 1992. Optimal strategies for a
generalized `Scissors, Paper and Stone' game.\ \textit{American Mathematical
Monthly} \textbf{99}: 935---942.

\bibitem {FisherRyan95a}Fisher, D., Ryan, J. 1995a. Tournament games and
positive tournaments.\ \textit{Journal of Graph Theory} \textbf{19}: 217---236.

\bibitem {FisherRyan95b}Fisher, D., Ryan, J. 1995b. Probabilities within
optimal strategies for tournament games.\ \textit{Discrete Applied
Mathematics} \textbf{56}: 87---91.

\bibitem {HofbauerSchlag00}Hofbauer, J., Schlag, K. 2000. Sophisticated imitation in cyclic games.\ \textit{Journal of Evolutionary Economics} \textbf{10}: 523---543.

\bibitem {LaffondetalGEB93}Laffond, G., Laslier, J.-F., Le Breton, M. 1993.
The Bipartisan set of a tournament game. \textit{Games and Economic Behavior}
\textbf{5}: 182---201.

\bibitem {Laslier97}Laslier, J.-F. 1997. \textit{Tournament Solutions and
Majority Voting}, Berlin: Springer-Verlag.

\bibitem {Laslier00a}Laslier, J.-F. 2000. Aggregation of preferences with a
variable set of alternatives. \textit{Social Choice and Welfare} \textbf{17}: 241---246.

\bibitem {Laslier00b}Laslier, J.-F. 2000. Interpretation of electoral mixed
strategies. \textit{Social Choice and Welfare} \textbf{17}: 247---267.

\bibitem {Levchenkov92}Levchenkov, V. S. 1992. Social choice theory: a new
insight. Discussion paper, Institute of Systems Analysis, Moscow.

\bibitem {McKelvey79}McKelvey, R. 1979. General conditions for global
intransitivities in a formal voting model.\ \textit{Econometrica} \textbf{47}: 1085---1112.

\bibitem {Moon68}Moon, J.W. 1968. \textit{Topics on Tournaments}. Holt,
Rinehart and Winston, New York.

\bibitem {Moulin86}Moulin, H. 1986. Choosing from a
tournament.\ \textit{Social Choice and Welfare} \textbf{3}: 271---291.

\bibitem {Myerson93}Myerson, R. B. 1993.\ Incentives to cultivate favored
minorities under alternative electoral systems. \textit{American Political
Science Review} \textbf{87}: 856---869.

\bibitem {Myerson95}Myerson, R. B. 1995.\ Analysis of democratic
institutions: structure, conduct and performance.\ \textit{Journal of Economic
Perspectives} \textbf{9}: 77---89.

\bibitem {Pemantle07}Pemantle, R. 2007.\ A survey of random processes with
reinforcement.\ \textit{Probability Surveys} \textbf{4}: 1---79.

\bibitem {Posch97}Posch, M. 1997.\ Cycling in a stochastic learning algorithm for normal form games.\ \textit{Journal of Evolutionary Economics} \textbf{7}: 193---207.

\bibitem {RivestChen10}Rivest, Ronald.L and Emily Shen (2010)
\textquotedblleft An Optimal Single-Winner Preferential Voting System Based on
Game Theory\textquotedblright\ http://people.csail.mit.edu/rivest/gt/latest\_conf.pdf

\bibitem {Rubinstein96}Rubinstein, Ariel 1996. Why are certain properties of
binary relations relatively more common in natural
languages?\ \textit{Econometrica} \textbf{64}: 343---355.

\bibitem {Shiryaev}Shiryaev, A.N. 1995. \textit{Probability}$.$\textit{
Graduate Text in Mathematics}, Springer.

\bibitem {SlutzkiVolij06}Slutzki, G. , Volij, O. 2006. Scoring of web pages
and tournaments --- axiomatizations.\ \textit{Social Choice and Welfare}
\textbf{26}: 75---92.

\bibitem {Ushakov76}Ushakov, I.\ A. 1976. The problem of choosing the
preferred element: An application to sport games. In \textit{Management
Science in Sports }(R. E. Machol, S. P. Ladany, and D.G.\ Morrison, eds.)
153---161.\ Amsterdam: North-Holland.
\end{thebibliography}
\end{document}